\documentstyle[11pt,amssymb]{article}
\newcommand{\lbrk}{{\linebreak[0]}}
\newcommand{\realtilde}{{\char'176}}
\newcommand{\ra}{\rightarrow}
\newcommand{\la}{\leftarrow}
\newcommand{\ts}{\textstyle}
\newcommand{\ds}{\displaystyle}

\newcommand{\artan}{{\rm artan}}
\newcommand{\W}{\hphantom{0}}

\newcommand{\slope}{
\setlength{\unitlength}{1mm}
\begin{picture}(3,12 )(0,0)
\put(0,-5){\line(1,4){2.75}}
\end{picture}}

\begin{document}

\vspace*{12mm}
\begin{center}
{\bf MULTIPLE-PRECISION ZERO-FINDING METHODS AND THE}\\
{\bf COMPLEXITY OF ELEMENTARY FUNCTION EVALUATION}%
\footnote{First appeared in
{\em Analytic Computational Complexity}
(edited by J F Traub), Academic Press, New York,
1975, 151--176.
Retyped with minor corrections and postscript
by Frances Page at Oxford University Computing
Laboratory, 1999 (urls updated 2005).\\
Copyright \copyright\ 1975--2010, R.~P.~Brent.
\hspace*{\fill} rpb028 typeset using \LaTeX.}\\[6ex]

{\bf Richard P Brent}\\[6ex]

{\bf Computer Centre,\\
 Australian National University,\\
Canberra, A.C.T.~2600, Australia}\\[6ex]

ABSTRACT
\end{center}

\begin{quote}
We consider methods for finding high-precision approximations to
simple zeros of smooth functions.  As an application, we give fast
methods for evaluating the elementary functions $\log(x), \exp(x),
\sin(x)$ etc.~to high precision.  For example, if $x$ is a positive
floating-point number with an $n$-bit fraction, then (under rather
weak assumptions) an $n$-bit approximation to $\log(x)$ or $\exp(x)$
may be computed 
in time 			%
asymptotically equal to $13M(n)\log_2n$ as $ n
\rightarrow \infty$, where $M(n)$ is the time required to multiply
floating-point numbers with $n$-bit fractions.  Similar results are
given for the other elementary functions, and some analogies with
operations on formal power series are mentioned.
\end{quote}

\section{Introduction}

When comparing methods for solving nonlinear equations or evaluating 
functions, it is customary to assume that the basic arithmetic 
operations (addition, multiplication, etc.) are performed with some
fixed precision.  However, an irrational number can only be 
approximated to arbitrary accuracy if the precision is allowed to
increase indefinitely.  Thus, we shall consider iterative processes
using variable precision.  Usually the precision will increase as the 
computation proceeds, and the final result will be obtained to high 
precision.  Of course, we could use the same (high) precision
throughout,   but then the computation would take longer than with
variable precision,  and the final result would be no more accurate.

\subsection*{Assumptions}

For simplicity we assume that a standard multiple-precision 
floating-point number representation is used, with a binary fraction of $n$
bits, where $n$ is large.  The exponent length is fixed, or may grow
as $o(n)$ if necessary.  To avoid table-lookup methods, we assume a 
machine with a finite random-access memory and a fixed number of 
sequential tape units.  Formally, the results hold for multitape
Turing machines.

\subsection*{Precision $n$ Operations}

An operation is performed with precision $n$ if the operands and
result are floating-point numbers as above (i.e., precision $n$
numbers), and the relative error in the result is $O(2^{-n})$.

\subsection*{Precision $n$ Multiplication}

Let $M(n)$ be the time required to perform precision $n$ multiplication.
(Time may be regarded as the number of single-precision operations,
or the number of bit operations, if desired.)  The classical method
gives $M(n) = O(n^2)$, but methods which are faster for large $n$ are
known.  Asymptotically the fastest method known is that of Sch\"{o}nhage
and Strassen $[71]$, which gives \vspace*{-4mm}

$$ M(n) = O(n\log(n)\log\log(n)) \;\;{\rm as}\;\; n \rightarrow \infty \;.
\leqno{(1.1)} $$

Our results do not depend on the algorithm used for multiplication,
provided $M(n)$ satisfies the following two conditions.\vspace*{-4mm}

$$ n = o(M(n))\;, \;\;{\rm i.e.,}\;\; \lim_{n \ra \infty} n/M(n) = 0 \;;
\leqno{(1.2)} $$

and, for any $\alpha > 0$,\vspace*{-4mm}

$$ M(\alpha n) \sim \alpha M(n)\;, \;\;{\rm i.e.,}\;\; 
            \lim_{n \ra \infty} \frac{M(\alpha n)}{\alpha M(n)} = 1 \;.
\leqno{(1.3)} $$

Condition (1.2) enables us to neglect additions, since the time for
an addition is $O(n)$, which is asymptotically negligible compared
to the time for a multiplication.  Condition (1.3) certainly holds
if\vspace*{-4mm}

$$ M(n) \sim cn[\log(n)]^\beta [\log\log(n)]^\gamma \;, \eqno{}$$

though it does not hold for some implementations of the 
Sch\"{o}nhage-Strassen method.  We need (1.3) to estimate the
constants in the asymptotic ``$O$'' results:  if the constants are
not required then much weaker assumptions suffice, as in Brent
$[$75a,b$]$.\\

The following lemma follows easily from (1.3).

\subsection*{Lemma 1.1}

If $0 < \alpha < 1, M(n) = 0$ for $n < 1$, and
$c_1 < \frac{1}{1-\alpha} < c_2$, then \vspace*{2mm}

$$ c_1M(n) < \sum^{\infty}_{k=0} M(\alpha^k n) < c_2M(n)$$

for all sufficiently large $n$.

\section{Basic Multiple-precision Operations}

In this section we summarize some results on the time required to
perform the multiple-precision operations of divison, extraction
of square roots, etc.  Additional results are given in Brent
$[$75a$]$.

\subsection*{Reciprocals}

Suppose $a \neq 0$ is given and we want to evaluate a precision
$n$ approximation to $1/a$.  Applying Newton's method to the equation
\vspace*{-4mm}

$$ f(x) \equiv a - 1/x = 0$$

gives the well-known iteration\vspace*{-4mm}

$$ x_{i+1} = x_i - x_i \varepsilon_i \;, $$

where\vspace*{-4mm}

$$ \varepsilon_i = ax_i - 1 \;. \eqno{}$$

Since the order of convergence is two, only $k \sim \log_2n$ iterations
are required if $x_0$ is a reasonable approximation to $1/a$,
e.g., a single-precision approximation.\\

If
$\varepsilon_{k-1}=O(2^{-n/2})$,	%
then
$\varepsilon_k  = O(2^{-n})$,
so at the last iteration it is sufficient
to perform the multiplication of $x_{k-1}$ by $\varepsilon_{k-1}$
using precision $n/2$, even though $ax_{k-1}$ must be evaluated with
precision $n$.  Thus, the time required for the last iteration is
$M(n) + M(n/2) + O(n)$.  The time for the next to last iteration is
$M(n/2) + M(n/4) + O(n/2)$, since this iteration need only give 
an approximation accurate to $O(2^{-n/2})$, and so on.  Thus,
using Lemma 1.1, the total time required is 

$$ I(n) \sim (1 + {\ts \frac{1}{2}})(1 + {\ts\frac{1}{2}} + {\ts\frac{1}{4}}
               + \cdots)M(n) \sim 3M(n)$$

as $n \ra \infty$.

\subsection*{Division}

Since $b/a = b(1/a)$, precision $n$ division may be done in time \vspace*{-2mm}

$$ D(n) \sim 4M(n) $$

as $n \ra \infty$.

\subsection*{Inverse Square Roots}

Asymptotically the fastest known method for evaluating  $a^{-\frac{1}{2}}$
to precision $n$ is to use the third-order iteration\vspace*{-4mm}

$$x_{i+1} = x_i - {\ts\frac{1}{2}}x_i(\varepsilon_i - {\ts\frac{3}{4}}
                          \varepsilon^2_i) \;, $$

where\vspace*{-4mm}

$$ \varepsilon_i = ax^2_i - 1 \;. \eqno{}$$

At the last iteration it is sufficient to evaluate $ax^2_i$ to
precision $n$, $\varepsilon^2_i$ to precision $n/3$, and
$x_i(\varepsilon_i - \frac{3}{4} \varepsilon^2_i)$ to precision
$2n/3$.  Thus, using Lemma 1.1 as above, the total time required is

$$ Q(n) \sim (2+{\ts\frac{1}{3}}+{\ts\frac{2}{3}})(1+{\ts\frac{1}{3}}
               +{\ts\frac{1}{9}}+\cdots) M(n) \sim {\ts\frac{9}{2}}M(n) $$

as $n \ra \infty$.

\subsection*{Square Roots}

Since \vspace*{-4mm}

$$ \sqrt{a} = \left\{ \begin{array}{lll}
                     a.a^{-\frac{1}{2}} &{\rm if}& a>0 \;,\\
                     0                 &{\rm if}& a = 0 \;,\\
                     \end{array} \right. $$

we can evaluate $\sqrt{a}$ to precision $n$ in time
\vspace*{-2mm}

$$ R(n) \sim {\ts\frac{11}{2}}M(n)
\eqno{}$$

as $n \ra \infty$.  Note the direct use of Newton's method in the form
\vspace*{-4mm}

$$ x_{i+1} = {\ts\frac{1}{2}}(x_i + a/x_i)
\leqno{(2.1)} $$

or\vspace*{-5mm}

$$ x_{i+1} = x_i + \left( \frac{a- x^2_i}{2x_i} \right)
\leqno{(2.2)} $$

\vspace*{2mm}
is asymptotically slower, requiring time $\sim 8M(n)$ or $\sim 6M(n)$
respectively.

\section{Variable-precision Zero-finding Methods}

Suppose $\zeta \neq 0$ is a simple zero of the nonlinear equation

$$ f(x) = 0 \;.
\eqno{}$$

Here, $f(x)$ is a sufficiently smooth function which can be evaluated 
near $\zeta$, with absolute error $O(2^{-n})$, in time $w(n)$.
We consider some methods for evaluating $\zeta$ to precision $n$.
Since we are interested in results for very large $n$, the time
required to obtain a good starting approximation is neglected.

\subsection*{Assumptions}

To obtain sharp results we need the following two assumptions, which
are similar to (1.2) and (1.3):\vspace*{-4mm}

$$ M(n) = o(w(n)) \;, \;\;{\rm i.e.,}\;\; \lim_{n \ra \infty} M(n)/w(n) = 0 \;;
\leqno{(3.1)}$$

and, for some $\alpha \geq 1$ and all $\beta >0$,\vspace*{-4mm}

$$ w(\beta n) \sim \beta^\alpha w(n)
\leqno{(3.2)}$$

as $n \ra \infty$.\\

From (3.1), the time required for a multiplication is negligible compared
to the time for a function evaluation, if $n$ is sufficiently large.
(3.2) implies (3.1) if $\alpha > 1$, and (3.2) certainly holds if,
for example, \vspace*{-4mm}

$$ w(n) \sim cn^\alpha [\log(n)]^\gamma[\log\log(n)]^\delta \;.
\eqno{}$$

The next lemma follows from our assumptions in much the same way
as Lemma 1.1.

\subsection*{Lemma 3.1}

If $0<\beta<1, w(n) = 0$ for $n<1$, and \vspace*{-2mm}

$$ c_1 < 1/(1 - \beta^\alpha) < c_2 \;,
\eqno{}$$

then \vspace*{-4mm}

$$ c_1w(n) < \sum^\infty_{k=0} w(\beta^kn) < c_2w(n)
\eqno{}$$

for all sufficiently large $n$.

\subsection*{A Discrete Newton Method}

To illustrate the ideas of variable-precision zero-finding methods, 
we describe a simple discrete Newton method.  More efficient methods
are described in the next three sections, and in Brent~$[$75a$]$.\\

Consider the iteration \vspace*{-4mm}

$$ x_{i+1} = x_i - f(x_i)/g_i \;,
\eqno{}$$

where $g_i$ is a one-sided difference approximation to $f'(x_i)$,
i.e., \vspace*{1mm}

$$ g_i = \frac{f(x_i + h_i) - f(x_i)}{h_i} \;.
\eqno{}$$

If $\varepsilon = \vert x_i - \zeta\vert$ is sufficiently small,
$f(x_i)$ is evaluated with absolute error $O(\varepsilon^2_i)$,
and $h_i$ is small enough that \vspace*{-4mm}

$$ g_i = f'(x_i) + O(\varepsilon_i) \;,
\leqno{(3.3)}$$

then the iteration converges to $\zeta$ with order at least two.
To ensure (3.3), take $h_i$ of order $\varepsilon_i$, e.g. $h_i=f(x_i)$.\\

To obtain $\zeta$ to precision $n$, we need two evaluations of $f$ 
with absolute error $O(2^{-n})$, preceded by two evaluations with
error, $O(2^{-n/2})$, etc.  Thus the time required is \vspace*{-4mm}

$$ t(n) \sim 2(1+2^{\alpha}+ 2^{-2\alpha}+ \cdots)w(n) \;.
\leqno{(3.4)}$$

\subsection*{Asymptotic Constants}

We say that a zero-finding method has {\em asymptotic constant\/} $C(\alpha)$
if, to find a simple zero $\zeta \neq 0$ to precision $n$, the
method requires time $t(n) \sim C(\alpha)w(n)$ as $n \ra \infty$.
(The asymptotic constant as defined here should not be confused with
the ``asymptotic error constant'' as usually defined for 
single-precision zero-finding methods.)\\

For example, from (3.4), the discrete Newton method described above
has asymptotic constant

$$ C_N(\alpha) = 2/(1 - 2^{-\alpha}) \leq 4 \;.
\eqno{}$$

Note that the time required to evaluate $\zeta$ to precision $n$
is only a small multiple of the time required to evaluate $f(x)$
with absolute error $O(2^{-n})$.  (If we used fixed precision, the
time to evaluate $\zeta$ would be of order $\log(n)$ times the time to
evaluate $f(x)$.)

\section{A Variable-precision Secant Method}

The secant method is known to be more efficient than the discrete Newton
method when fixed-precision arithmetic is used.  The same is true with
variable-precision arithmetic, although the ratio of efficiencies is no
longer constant, but depends on the exponent $\alpha$ in (3.2).  Several
secant-like methods are described in Brent $[$75a$]$; here we consider the
simplest such method, which is also the most efficient if $\alpha <
4.5243 \ldots$\\

The secant iteration is \vspace*{-3mm}

$$ x_{i+1} = x_i - f_i\left(\frac{x_i-x_{i-1}}{f_i - f_{i-1}} \right) \;,
\eqno{}$$

\vspace*{1mm}
where $f_i = f(x_i)$, and we assume that the function evaluations are
performed with suficient accuracy to ensure that the order of convergence
is at least $\rho = (1+\sqrt{5})/2 = 1.6180 \ldots ,$ the
larger root of \vspace*{-4mm}

$$ \rho^2 = \rho + 1 \;.
\leqno{(4.1)} $$

Let $\varepsilon = \vert x_{i-1} - \zeta \vert$.  Since the smaller root
of (4.1) lies inside the unit circle, we have

$$ x_i - \zeta = O(\varepsilon^\rho)
\eqno{}$$

\vspace*{-2mm}
and \vspace*{-5mm}

$$ x_{i+1} - \zeta = O(\varepsilon^{\rho^2}) \;.
\eqno{} $$

To give order $\rho$, $f_i$ must be evaluated with absolute error
$O(\varepsilon^{\rho^2})$.  Since $f_i = O(\vert x_i - \zeta\vert) =
O(\varepsilon^\rho)$, it is also necessary to evaluate
$(f_i - f_{i-1})/(x_i - x_{i-1})$ with relative error
$O(\varepsilon^{\rho^2 - \rho})$, but $\vert x_i - x_{i-1}\vert \sim
\varepsilon$, so it is necessary to evaluate $f_{i-1}$ with absolute
error $O(\varepsilon^{\rho^2 - \rho + 1})$.  $[$ Since $f_i$ must be
evaluated with absolute error $O(\varepsilon^{\rho^2})$, $f_{i-1}$
must be evaluated with absolute error $O(\varepsilon^\rho)$, but
$\rho^2 - \rho + 1 = 2 > \rho$, so this condition is superfluous.$]$\\

The conditions mentioned are sufficient to ensure that the order of
convergence is at least $\rho$.  Thus, if we replace 
$\varepsilon^{\rho^2}$ by $2^{-n}$, we see that $\zeta$ may be
evaluated to precision $n$ if $f$ is evaluated with absolute errors
$O(2^{-n})$, $O(2^{-2n\rho^{-2}})$, $O(2^{-2n\rho^{-3}})$,
$O(2^{-2n\rho^{-4}})$, \ldots $\;\;$It follows that the asymptotic constant
of the secant method is\vspace*{-4mm}

$$ C_S(\alpha) = 1 + (2\rho^{-2})^\alpha/(1-\rho^{-\alpha}) \leq
                 C_S(1) = 3 \;.
\eqno{}$$

The following lemma states that the secant method is asymptotically
more efficient than the discrete Newton method when variable precision
is used.

\subsection*{Lemma 4.1}

$C_S(\alpha) < C_N(\alpha)$ for all $\alpha \geq 1$.  In fact
$C_S(\alpha)/C_N(\alpha)$ decreases monotonically from $\frac{3}{4}$
(when $\alpha = 1$) to $\frac{1}{2}$ (as $\alpha \ra \infty$).

\section{Other Variable-precision Interpolatory Methods}

With fixed precision, inverse quadratic interpolation is more efficient
than linear interpolation, and inverse cubic interpolation is even more
efficient, if the combinatory cost (i.e., ``overhead'') is negligible.
With variable precision the situation is different.  Inverse quadratic
interpolation is slightly more efficient than the secant method, but
inverse cubic interpolation is {\em not\/} more efficient than inverse
quadratic interpolation if $\alpha \leq 4.6056\ldots$  $\;\;$  Since  the
combinatory cost of inverse cubic interpolation is considerably higher
than that of inverse quadratic interpolation, the inverse cubic method 
appears even worse if combinatory costs are significant.

\subsection*{Inverse Quadratic Interpolation}

The analysis of variable-precision methods using inverse quadratic
interpolation is similar to that for the secant method, so we only
state the results.  The order $\rho = 1.8392\ldots$ $\;\;$ is the positive 
root of $\rho^3 = \rho^2+\rho+1$.  It is convenient to define
$\sigma  = 1/\rho   = 0.5436\ldots$ $\;\;$ To evaluate $\zeta$ to precision
$n$ requires evaluations of $f$ to (absolute) precision $n$,
$(1-\sigma +\sigma^2)n$, and $\sigma^j(1 - \sigma -\sigma^2 + 2\sigma^3)n$
for $j=0,1,2,\ldots$ $\;\;$ Thus, the asymptotic constant is 

$$ \begin{array}{lll}

C_Q(\alpha)& =& 1+(1-\sigma+\sigma^2)^\alpha + (3\sigma^3)^\alpha/
                                   (1-\sigma^\alpha)\\[1ex]
           &\leq& C_Q(1) = {\ts\frac{1}{2}}(7-2\sigma-\sigma^2) = 2.8085\ldots  \;.
   \end{array}
\eqno{}$$

\subsection*{Lemma 5.1}

$C_Q(\alpha) < C_S(\alpha)$ for all $\alpha \geq 1$.  In fact,
$C_Q(\alpha)/C_S(\alpha)$ increases monotonically from 0.9361\ldots $\;\;$
(when $\alpha = 1$) to 1 (as $\alpha \ra \infty$).

\subsection*{Inverse Cubic Interpolation, etc}
   
If $\mu = 0.5187\ldots$ $\;\;$is the positive root of $\mu^4+\mu^3+\mu^2+\mu
= 1$, then the variable-precision  method of order $1/\mu = 1.9275\ldots \;,$
using inverse cubic interpolation, has asymptotic constant

$$ \begin{array}{lll}

C_C(\alpha) &=& 1+(1-\mu + \mu^2)^\alpha + 
  (1 - \mu - \mu^2 +2\mu^3)^\alpha +
  (4\mu^4)^\alpha/(1-\mu^\alpha)\\[1ex]

       &\leq& C_C(1) = (13 - 6\mu - 4\mu^2 - 2\mu^3)/3 = 2.8438\ldots \;.
\end{array}
\eqno{}$$

\vspace*{2mm}
Note that $C_C(1) > C_Q(1)$.  Variable-precision methods using inverse
interpolation of arbitrary degree are described in Brent $[$75a$]$.
Some of these methods are slightly more efficient than the inverse 
quadratic interpolation method if $\alpha$ is large, but inverse
quadratic interpolation is the most efficient method known for
$\alpha < 4.6056\ldots$ $\;\;$ In practice $\alpha$ is usually 1,
$1\frac{1}{2}$ or 2.

\subsection*{An Open Question}

Is there a method with asymptotic constant $C(\alpha)$ such that
$C(1) < C_Q(1)$?

\section{Variable-precision Methods using Derivatives}

In Sections 3 to 5 we considered methods for solving the nonlinear
equation $f(x)=0$, using only evaluations of $f$.  Sometimes it is 
easy to evaluate $f'(x), f''(x), \ldots$ $\;$ once $f(x)$ has been 
evaluated, and the following theorem shows that it is possible to
take advantage of this.  For an application, see Section~10.

\subsection*{Theorem 6.1}

If the time to evaluate $f(x)$ with an absolute error $O(2^{-n})$
is $w(n)$, where $w(n)$ satisfies conditions (3.1) and (3.2), and
(for $k=1,2,\ldots$) the time to evaluate $F^{(k)}(x)$ with absolute
error $O(2^{-n})$ is $w_k(n)$, where \vspace*{-4mm}

$$ w_k(n) = o(w(n))
\eqno{}$$

as $n \ra \infty$, then the time to evaluate a simple zero $\zeta \neq
0$ of $f(x)$ to precision $n$ is

$$ t(n) \sim w(n)
\eqno{}$$

\vspace*{-2mm}
as $n \ra \infty$.

\subsubsection*{Proof}

For fixed $k\geq1$, we may use a direct or inverse Taylor series
method of order $k+1$.  The combinatory cost is of order
$k\log(k+1)M(n)$ (see Brent and Kung $[$75$]$).  From (3.1),
this is $o(w(n))$ as $n \ra \infty$.  Thus, \vspace*{-5mm}

\begin{eqnarray}
  t(n) &\leq& [1-(k+1)^{-\alpha}]^{-1}w(n)+o(w(n))\nonumber \\ [1ex]
       &\leq& (1+{\textstyle\frac{1}{k}}+o(1))w(n) \;. \nonumber
\end{eqnarray}\\[-5ex]

For sufficiently large $n$, the ``$o(1)$'' term is less than $1/k$,
so

$$ t(n) \leq (1+\ts\frac{2}{k})w(n) \;.
\eqno{}$$

Given $\varepsilon > 0$, choose $k \geq 2/\varepsilon$.  Then, for
all sufficiently large $n$,

$$ w(n) \leq t(n) \leq (1+\varepsilon)w(n) \;,
\eqno{}$$

\vspace*{-2mm}
so $t(n) \sim w(n)$ as $n \ra \infty$.

\subsection*{Corollary 6.1}

If the conditions of Theorem 6.1 hold, $f\colon[a,b] \ra I$, $f'(x) \neq
0$ for $x \in [a,b]$, and $g$ is the inverse function of $f$, then
the time to evaluate $g(y)$ with absolute error $O(2^{-n})$, for
$y\in I$, is

$$ w_g(n) \sim w(n)
\eqno{}$$

\vspace*{-2mm}
as $n \ra \infty$.

\subsubsection*{Note}

Corollary 6.1 is optimal in the sense that, if $w_g(n) \sim cw(n)$
for some constant $c<1$, then $w(n) \sim cw_g(n)$ by the same argument,
so $w(n) \sim c^2w(n)$, a contradiction.  Hence, $c=1$ is minimal.

\section{The Arithmetic-geometric Mean Iteration}

Before considering the multiple-precision evaluation of elementary
functions, we recall some properties of the arithmetic-geometric 
(A--G) mean iteration of Gauss $[$1876$]$.  Starting from any two
positive numbers $a_0$ and $b_0$, we may iterate as follows:

$$ a_{i+1} = \frac{a_i+b_i}{2} \hspace*{8mm}\makebox{arithmetic mean}
\eqno{}$$

\vspace*{-3mm}
and
\vspace*{-3mm}

$$ b_{i+1}= \sqrt{a_ib_i} \hspace*{8mm}\makebox{geometric mean}
\eqno{}$$

\vspace*{-3mm}
for $i=0,1,\ldots$

\subsection*{Second-order Convergence}

The A--G mean iteration is of computational interest because it 
converges very fast.  If $b_i \ll a_i$, then

$$ b_{i+1}/a_{i+1} = \frac{2\sqrt{b_i/a_i}}
                          {1+ b_i/a_i}
         \simeq 2\sqrt{b_i/a_i}  \;,
\eqno{}$$

\vspace*{1mm}
so only about $\vert \log_2(a_0/b_0)\vert$ iterations are required before
$a_i/b_i$ is of order 1.  Once $a_i$ and $b_i$ are close together the
convergence is second order, for if $b_i/a_i = 1-\varepsilon_i$ then
$$\varepsilon_{i+1}= 1-b_{i+1}/a_{i+1}=1-2(1-\varepsilon_i)^\frac{1}{2}/
          (2-\varepsilon_i)=\varepsilon^2_i/8+O(\varepsilon^3_i)\;.$$

\subsection*{Limit of the A--G Mean Iteration}

There is no essential loss of generality in assuming that $a_0=1$
and $b_0=\cos\phi$ for some $\phi$.  If $a={\ds\lim_{1\ra \infty}} a_i =
{\ds\lim_{i\ra \infty}} b_i$, then \vspace*{-2mm}

$$ a  = \frac{\pi}{2K(\phi)} \;,
\leqno{(7.1)}$$

where $K(\phi)$ is the complete elliptic integral of the first kind,
i.e.,

$$ K(\phi) = \int\limits^{\pi/2}_0
		\frac{1}{\sqrt{1- \sin^2\phi\sin^2\theta}}\;
                d\theta \;.
\eqno{}$$

\vspace*{2mm}
(A simple proof of (7.1) is given in Melzak $[$73$]$.)\\

Also, if $c_0=\sin\phi$, $c_{i+1}=a_i - a_{i+1}$ $(i=0,1,\ldots)$,
then \vspace*{-3mm}

$$ \sum^\infty_{i=0} 2^{i-1}c^2_i = 1 - \frac{E(\phi)}{K(\phi)} \;,
\leqno{(7.2)}$$

\vspace*{1mm}
where $E(\phi)$ is the complete elliptic integral of the second kind, i.e.,

$$ E(\phi) = \int\limits^{\pi/2}_0
		\sqrt{1- \sin^2\phi\sin^2\theta}\; d\theta \;.
\eqno{}$$

\vspace*{2mm}
Both (7.1) and (7.2) were known by Gauss.

\subsection*{Legendre's Identity}

For future use, we note the identity \vspace*{-4mm}

$$ K(\phi)E(\phi') + K(\phi')E(\phi) - K(\phi)K(\phi') = \ts\frac{1}{2}\pi \;,
\leqno{(7.3)}$$

where $\phi + \phi' = \frac{1}{2}\pi$.  (Legendre $[$11$]$ proved by
differentiation that the left side of (7.3) is constant, and the
constant may be determined by letting $\phi \ra 0$.)

\section{Fast Multiple-precision Evaluation of $\pi$}
\label{sec:pi}

The classical methods for evaluating $\pi$ to precision $n$ take time
$O(n^2)$:  see, for example, Shanks and Wrench $[$62$]$.  Several 
methods which are asymptotically faster than $O(n^2)$ are known.  For 
example, in Brent $[$75a$]$ a method which requires  time $O(M(n)\log^2(n))$
is described.  From the bound (1.1) on $M(n)$, this is faster than
$O(n^{1+\varepsilon})$ for any $\varepsilon > 0$.\\

Asymptotically the fastest known methods require time $O(M(n)\log(n))$.  
One such method is sketched in Beeler et al $[$72$]$.  The method given 
here is faster, and does not require the preliminary computation of $e$.

\subsection*{The Gauss-Legendre Method}

Taking $\phi = \phi' = \pi/4$ in (7.3), and dividing both sides by
$\pi^2$, we obtain \vspace*{-4mm}

$$ [2K(\pi/4)E(\pi/4)-K^2(\pi/4)]/\pi^2 = \frac{1}{2\pi} \;.
\leqno{(8.1)}$$

However, from the A--G mean iteration with $a_0 =1$ and 
$b_0=2^{-\frac{1}{2}}$, and the relations (7.1) and (7.2), we can evaluate
$K(\pi/4)/\pi$ and $E(\pi/4)/\pi$, and thus the left side of (8.1).
A division then gives $\pi$.  (The idea of using (7.3) in this way
occurred independently to Salamin $[$75$]$ and Brent $[$75b$]$.)
After a little simplification, we obtain the following algorithm (written
in pseudo-Algol):\\

\hspace*{12mm}$A\la 1;\; B\la 2^{-\frac{1}{2}};\; T \la 1/4;\; X \la 1;$\\

\hspace*{12mm}$\begin{array}{ll}
\mbox{\bf while} & A - B > 2^{-n} \;\; \mbox{\bf do}\\
                 & \mbox{\bf begin} \;\; Y \la A;\; A \la \frac{1}{2}(A+B);\;
                            B \la \sqrt{BY} \;;\\
                 & \hspace*{12mm} T \la T - X(A - Y)^2 \;;\\
                 & \hspace*{12mm} X \la 2X\\
                 & \mbox{\bf end};\\
\mbox{\bf return}& A^2/T \;\; [\mbox{or, better,}\; (A+B)^2/(4T)] \;.
\end{array}$\\\\

The rate of convergence is illustrated in Table 8.1.

\subsection*{Table 8.1:  Convergence of the Gauss-Legendre Method}

\vspace*{3mm}
\begin{center}
\begin{tabular}{cclcl}

Iteration && $A^2/T - \pi$ && \multicolumn{1}{c}{$\pi - (A+B)^2/(4T)$} \\ \hline

0&&8.6$ \; \times \; 10^{-1}$  &&2.3$ \; \times \; 10^{-1}$\\
1&&4.6$ \; \times \; 10^{-2}$  &&1.0$ \; \times \; 10^{-3}$\\
2&&8.8$ \; \times \; 10^{-5}$  &&7.4$ \; \times \; 10^{-9}$\\
3&&3.1$ \; \times \; 10^{-10}$ &&1.8$ \; \times \; 10^{-19}$\\
4&&3.7$ \; \times \; 10^{-21}$ &&5.5$ \; \times \; 10^{-41}$\\
5&&5.5$ \; \times \; 10^{-43}$ &&2.4$ \; \times \; 10^{-84}$ \\
6&&1.2$ \; \times \; 10^{-86}$ &&2.3$ \; \times \; 10^{-171}$\\
7&&5.8$ \; \times \; 10^{-174}$&&1.1$ \; \times \; 10^{-345}$\\
8&&1.3$ \; \times \; 10^{-348}$&&1.1$ \; \times \; 10^{-694}$\\
9&&6.9$ \; \times \; 10^{-698}$&&6.1$ \; \times \; 10^{-1393}$\\\hline\hline
\end{tabular}
\end{center}

\vspace*{4mm}
Since the A--G mean iteration converges with order 2, we need
$\sim \log_2n$ iterations to obtain precision $n$.  Each iteration
involves one (precision $n$) square root, one multiplication, one
squaring, one multiplication by a power of two, and some additions.
Thus from the results of Section~2, the time required to evaluate
$\pi$ is $\sim \frac{15}{2}M(n)\log_2n$.

\subsection*{Comments}

\begin{enumerate}
\item Unlike Newton's iteration, the A--G mean iteration is not
self-correcting.  Thus, we cannot start with low precision and
increase it, as was possible in Section~2.
\item Since there are $\sim\log_2n$ iterations, we may lose
$O(\log\log(n))$ bits of accuracy through accumulation of rounding
errors, even though the algorithm is numerically stable.  Thus,
it may be necessary to work with precision $n+O(\log\log(n))$.
From (1.3), the time required is still $\sim\frac{15}{2}M(n)\log_2n$.
\end{enumerate}

\section{Multiple-precision Evaluation of $\log(x)$}

There are several algorithms for evaluating $\log(x)$ to precision $n$
in time $O(M(n)\log(n))$.  For example, a method based on Landen
transformations of incomplete elliptic integrals is described in 
Brent $[$75b$]$.  The method described here is essentially due to
Salamin (see Beeler et al $[$72$]$), though the basic relation (9.1)
was known by Gauss.\\

If $\cos(\phi) = \varepsilon^\frac{1}{2}$ is small, then \vspace*{-4mm}

$$ K(\phi) = (1 + O(\varepsilon))\log(4\varepsilon^{-\frac{1}{2}})
\leqno{(9.1)}$$

Thus, taking $a_0=1, b_0=4/y$, where $y=4\varepsilon^{-\frac{1}{2}}$,
and applying the A--G mean iteration to compute $a={\ds\lim_{i\ra \infty}}a_i$,
gives \vspace*{-3mm}

$$ \log(y) = \frac{\pi}{2a}(1+O(y^{-2}))
\eqno{}$$

for large $y$.  Thus, so long as $y\geq 2^{n/2}$, we can evaluate
$\log(y)$ to precision $n$.  If $\log(y)=O(n)$ then $\sim2\log_2n$
iterations are required, so the time is $\sim13M(n)\log_2n$, assuming
$\pi$ is precomputed.\\

For example, to find $\log(10^6)$ we start the A--G mean iteration with
$a_0=1$ and $b_0=4 \times10^{-6}$.  Results of the first seven iterations are given
to 10 significant figures in Table 9.1.  We find that $\pi/(2a_7)
= 13.81551056$, which is correct.

\subsection*{Table 9.1:  Computation of $\log(10^6)$}

\vspace*{3mm}
\begin{center}
\begin{tabular}{ccc}

$i$ & $a_i$ & $b_i$ \\ \hline

0 &$\,\,\W1.000000000 \;\;\hphantom{\times 10^{-000}} $
	 		       &$4.000000000 \times 10^{-6}$\\
1 &$5.000020000 \times 10^{-1}$&$2.000000000 \times 10^{-3}$\\
2 &$2.510010000 \times 10^{-1}$&$3.162283985 \times 10^{-2}$\\
3 &$1.413119199 \times 10^{-1}$&$8.909188753 \times 10^{-2}$\\
4 &$1.152019037 \times 10^{-1}$&$1.122040359 \times 10^{-1}$\\
5 &$1.137029698 \times 10^{-1}$&$1.136930893 \times 10^{-1}$\\
6 &$1.136980295 \times 10^{-1}$&$1.136980294 \times 10^{-1}$\\
7 &$1.136980295 \times 10^{-1}$&$1.136980295 \times 10^{-1}$\\\hline\hline
\end{tabular}
\end{center}

\vspace*{4mm}

Since $\log(2) = \frac{1}{n}\log(2^n)$, we can evaluate $\log(2)$
to precision $n$ in time $\sim13M(n)\log_2n$.  Suppose $ x \in [b,c]$,
where $b>1$.  We may set $y=2^nx$, evaluate $\log(y)$ as above, and
use the identity

$$\log(x) = \log(y) - n\log(2)
\eqno{}$$

to evaluate $\log(x)$.  Since $\log(y) \simeq n\log(2)$, approximately
$\log_2n$ significant bits will be lost through cancellation, so it is
necessary to work with precision $n+O(\log(n))$.\\

If $x$ is very close to $1$, we have to be careful in order to obtain
$\log(x)$ with a small {\em relative\/} error.  Suppose
$x=1+\delta$.  If $\vert\delta\vert < 2^{-n/\log(n)}$ we may use the
power series

$$ \log(1+\delta) = \delta - \delta^2/2 + \delta^3/3 - \, \ldots \;,
\eqno{}$$

and it is sufficient to take about $\log(n)$ terms.  If $\delta$
is larger, we may use the above A--G mean method, with working
precision $n+O(n/\log(n))$ to compensate for any cancellation.\\

Finally, if $0<x<1$, we may use $\log(x)= -\log(1/x)$, where
$\log(1/x)$ is computed as above.  To summarize, we have proved:

\subsection*{Theorem 9.1}

If $x>0$ is a precision $n$ number, then $\log(x)$ may be evaluated
to precision $n$ in time\linebreak
$\sim13M(n)\log_2n$ as $n\ra \infty$
$[$ assuming $\pi$ and $\log(2)$ precomputed to precision 
$n+O(n/\log(n))$$]$.

\subsection*{Note}

The time required to compute $\log(x)$ by  the obvious power series
method is $O(nM(n))$.  Since $13\log_2n < n$  for $n\geq 83$, the
method described here may be useful for moderate $n$, even if the
classical $O(n^2)$ multiplication algorithm is used.

\section{Multiple-precision Evaluation of $\exp(x)$}

Corresponding to Theorem 9.1 we have:

\subsection*{Theorem 10.1}

If $[a,b]$ is a fixed interval, and $x \in [a,b]$ is a precision $n$
number such that $\exp(x)$ does not underflow or overflow, the
$\exp(x)$ can be evaluated to precision $n$ in time $\sim13M(n)\log_2n$
as $n\ra \infty$ (assuming $\pi$ and $\log(2)$ are precomputed).

\subsubsection*{Proof}

To evaluate $\exp(x)$ we need to solve the equation $f(y)=0$, where
$f(y)=\log(y)-x$, and $x$ is regarded as constant.  Since

$$
f^{(k)}(y)=(-1)^{k-1} \; (k-1)\,!\,y^{-k} \eqno{}
$$

can be evaluated in time $O(M(n)) = o(M(n)\log(n))$ for any fixed
$k \geq 1$, the result follows from Theorems 6.1 and 9.1.
We remark that the $(k+1)$-th order method in the proof of Theorem 6.1 may
simply be taken as \vspace*{-3mm}

$$ y_{i+1}= y_i \sum^k_{j=0}(x-\log(y_i))^j/j\,!\,
\eqno{}$$

\section{Multiple-precision Operations on Complex Numbers}

Before considering the multiple-precision evaluation of trigonometric
functions, we need to state some results on multiple-precision operations
with complex numbers.  We assume that a precision $n$ complex number
$z=x+iy$ is represented as a pair $(x, y)$ of precision $n$ real
numbers.  As before, a precision $n$ operation is one which gives a
result with a relative error $O(2^{-n})$.  (Now, of course, the relative
error may be complex, but its absolute value must be $O(2^{-n})$.)
Note that the smaller component of a complex result may occasionally
have a large relative error, or even the wrong sign!

\subsection*{Complex Multiplication}

Since $z=(t+iu)(v+iw)=(tv-uw)+i(tw+uv)$, a complex multiplication may be 
done with four real multiplications and two additions.  However, we may 
use an idea of Karatsuba and Ofman $[$62$]$ to reduce the work required
to three real multiplications and some additions:  evalute $tv, uw$,
and $(t+u)(v+w)$, then use \vspace*{-3mm}

$$ tw+uv=(t+u)(v+w)-(tv+uw) \;.
\eqno{}$$

Since $\vert t+u\vert \;\leq \sqrt{2}\;\vert t+iu\vert$ and
$\vert v+w\vert \; \leq \sqrt{2} \;\vert v+iw\vert$, we have

$$\vert(t+u)(v+w)\vert \;\leq\; 2\vert z\vert \;.
\eqno{}$$

Thus, all rounding errors are of order $2^{-n}\vert z\vert$ or less,
and the computed product has a relative error $O(2^{-n})$.  The time
for the six additions is asymptotically negligible compared to that
for the three multiplications, so precision $n$ complex multiplication
may be performed in time $\sim3M(n)$.

\subsection*{Complex Squares}

Since $(v+iw)^2 = (v-w)(v+w)+2ivw$, a complex square may be evaluated with
two real multiplications and additions, in time $\sim2M(n)$.

\subsection*{Complex Division}

Using complex multiplication as above, and the same division algorithm
as in the real case, we can perform complex division in time $\sim12M(n)$.
However, it is faster to use the identity

$$\frac{t+iu}{v+iw}=(v^2+w^2)^{-1}[(t+iu)(v-iw)] \;,
\eqno{}$$

reducing the problem to one complex multiplication, four real 
multiplications, one real reciprocal, and some additions.  This gives time 
$\sim10M(n)$.  For complex reciprocals we have $t=1, u=0$, and time
$\sim7M(n)$.

\subsection*{Complex Square Roots}

Using (2.2) requires, at the last iteration, one precision $n$ complex
squaring and one precision $n/2$ complex division.  Thus, the time
required is $\sim 2(2+ 10/2)M(n) = 14M(n)$.

\subsection*{Complex A--G Mean Iteration}

From the above results, a complex square root and multiplication may be
performed in time $\sim 17M(n)$.  Each iteration transforms two points 
in the complex plane into two new points, and has an interesting
geometric interpretation.

\section{Multiple-precision Evaluation of Trigonometric Functions}

Since \vspace*{-4mm}

$$ \log(v+iw)=\log\vert v+iw\vert \;+ i.\artan(w/v)
\leqno{(12.1)}$$

and \vspace*{-4mm}

$$ \exp(i\theta) = \cos(\theta) + i.\sin(\theta) \;,
\leqno{(12.2)}$$

we can evaluate $\artan$$(x)$, $\cos(x)$ and $\sin(x)$ if we can evaluate
$\log(z)$ and $\exp(z)$ for complex arguments $z$.  This may be done 
just as described above for real $z$, provide we choose the correct
value of $\sqrt{a_j b_j}$.  Some care is necessary to avoid
excessive cancellation; for example, we should use the power series
for $\sin(x)$ if $\vert x\vert$ is very small, as described above for
$\log(1+\delta)$.  Since $\sim2\log_2n\;$ A--G mean iterations
are required to evaluate $\log(z)$, and each iteration requires time
$\sim17M(n)$, we can evaluate $\log(z)$ in time $\sim34M(n)\log_2n$.
From the complex version of Theorem 6.1, $\exp(z)$ may also be
evaluated in time $\sim34M(n)\log_2n$. \\

As an example, consider the evaluation of $\log(z)$ for $z=10^6(2+i)$.
The A--G mean iteration is started with $a_0=1$ and 
$b_0=4/z=1.6\times10^{-6} - (8.0\times10^{-7})i$.  The results of six iterations are given,
to 8 significant figures, in Table 12.1.

\subsection*{Table 12.1:  Evaluation of $\log10^6(2+i)$}

\vspace*{3mm}

\begin{tabular}{lllll}

$j$&&$a_j$&&$b_j$ \\[-2ex]
&&&&\\\hline
&&&&\\[-2ex]
0&$(1.0000000,$&$                0.0000000)$
 &$(1.6000000\times10^{-6},$&$\hspace*{-3mm}  -8.0000000\times10^{-7})$\\

1\hspace*{4mm}&$(5.0000080\times10^{-1},$&\hspace*{-3mm}$ -4.0000000\times10^{-7})$ \hspace*{4mm}
 &$(1.3017017\times10^{-3},$&\hspace*{-3mm}$ -3.0729008\times10^{-4})$\\

2&$(2.5065125\times10^{-1},$&\hspace*{-3mm}$ -1.5384504\times10^{-4})$
 &$(2.5686505\times10^{-2},$&\hspace*{-3mm}$ -2.9907884\times10^{-3})$\\

3&$(1.3816888\times10^{-1},$&\hspace*{-3mm}$ -1.5723167\times10^{-3})$
 &$(8.0373334\times10^{-2},$&\hspace*{-3mm}$ -4.6881008\times10^{-3})$\\

4&$(1.0927111\times10^{-1},$&\hspace*{-3mm}$ -3.1302088\times10^{-3})$
 &$(1.0540970\times10^{-1},$&\hspace*{-3mm}$ -3.6719673\times10^{-3})$\\

5&$(1.0734040\times10^{-1},$&\hspace*{-3mm}$ -3.4010880\times10^{-3})$
 &$(1.0732355\times10^{-1},$&\hspace*{-3mm}$ -3.4064951\times10^{-3})$\\

6&$(1.0733198\times10^{-1},$&\hspace*{-3mm}$ -3.4037916\times10^{-3})$
 &$(1.0733198\times10^{-1},$&\hspace*{-3mm}$ -3.4037918\times10^{-3})$\\ [-2ex]
&&&&\\\hline\hline
\end{tabular}\\\\

We find that \vspace*{-8mm}

\begin{eqnarray*}
\frac{\pi}{2a_7}&=& 14.620230 + 0.46364761i \\
               &\simeq& \log\vert z\vert + i.\artan({\ts\frac{1}{2}})
\end{eqnarray*}

as expected.\\

Another method for evaluating trigonometric functions in time
$O(M(n)\log(n))$, without using the identities (12.1) and (12.2)
is described in Brent $[$75b$]$.

\section{Operations on Formal Power Series}

There is an obvious similarity between a multiple-precision number
with base $\beta$: %

$$ \beta^e \sum^n_{i=1} a_i\beta^{-i}\;\; (0\leq a_i < \beta) \;,
\eqno{}$$

\vspace*{1mm}
and a formal power series:

$$ \sum^\infty_{i=0} a_ix^i \hspace*{12mm}
                  \mbox{($a_i$ real, $x$ an indeterminate)} \;.
\eqno{}$$

\vspace*{2mm}
Thus, it is not surprising that algorithms similar to those described
in Section~2 may be used to perform operations on power series.\\

In this section only, $M(n)$ denotes the number of scalar operations 
required to evaluate the first $n$ coefficients $c_0, \ldots , c_{n-1}$
in  the formal product \vspace*{1mm}

$$\left(\sum^\infty_{i=0}a_i x^i\right) 
   \left(\sum^\infty_{i=0} b_i x^i\right) =
   \sum^\infty_{i=0} c_ix^i \;.
\eqno{}$$

\vspace*{2mm}
Clearly, $c_j$ depends only on $a_0, \ldots , a_j$ and $b_0, \ldots , b_j$,
in fact

$$c_j = \sum^j_{i=0}a_i b_{j-i} \;.
\eqno{}$$

\vspace*{2mm}
The classical algorithm gives $M(n) = O(n^2)$, but it is possible to
use the fast Fourier transform (FFT) to obtain \vspace{-3mm}

$$ M(n) = O(n\log(n))
\eqno{}$$

\vspace{-2mm}
(see Borodin $[$73$]$).\\

If we assume that $M(n)$ satisfies conditions (1.2) and (1.3), then the
time bounds given in Section~2 for division, square roots, etc.~of
multiple-precision numbers also apply for the corresponding operations
on power series (where we want the first $n$ terms in the result).
For example, if 
$P(x) = {\ds \sum^\infty_{i=0}} a_i x^i$ and $a_0 \neq 0$, then the
first $n$ terms in the expansion of $1/P(x)$ may be found with
$\sim3M(n)$ operations as $n \ra \infty$.  However, some operations,
e.g.~computing exponentials, are much easier for power series than for
multiple-precision numbers!

\subsection*{Evaluation of $\log(P(x))$}

If $a_0 > 0$ we may want to compute the first $n$ terms in the power series
$Q(x) = \log(P(x))$.  Since $Q(x) = \log(a_0) + \log(P(x)/a_0)$, there is
no loss of generality in assuming that $a_0 = 1$.  Suppose
$Q(x)= {\ds \sum^\infty_{i=0}} b_i x^i$.  From the relation \vspace*{-1mm}

$$ Q'(x) = P'(x)/P(x) \;,
\leqno{(13.1)}$$

where the prime denotes formal differentiation with respect to $x$,
we have \vspace*{-8mm}

$$ \sum^\infty_{i=1} ib_i x^{i-1} =
\left( \sum^\infty_{i=1} ia_i x^{i-1}\right) \slope
\left( \sum^\infty_{i=0} a_i x^i\right) \;.
\leqno{(13.2)}$$

\vspace*{2mm}
The first $n$ terms in the power series for the right side of (13.2)
may be evaluated with $\sim 4M(n)$ operations, and then we need only
compare coefficients to find $b, \ldots , b_{n-1}$.  (Since
$a_0 = 1$, we know that $b_0 = 0$.)  Thus, the first $n$ terms in
$\log(P(x))$ may be found in $\sim 4M(n)$ operations.  It is
interesting to compare this result with Theorem 9.1.

\subsection*{Evaluation of $\exp(P(x))$}

If $R(x) = \exp(P(x))$ then $R(x) = \exp(a_0)\exp(P(x) - a_0)$,
so there is no loss of generality in assuming that $a_0 = 0$.
Now $\log(R(x)) - P(x) = 0$, and we may regard this as an equation 
for the unknown power series $R(x)$, and solve it by  one of the
usual iterative methods.  For example, Newton's method gives the
iteration \vspace*{-4mm}

$$ R_{i+1}(x) = R_i(x) - R_i(x)(\log(R_i(x)) - P(x)) \;.
\leqno{(13.3)}$$

If we use the starting approximation $R_0(x)=1$, then the terms in
$R_k(x)$ agree exactly with those in $R(x)$ up to (but excluding)
the term $O(x^{2^k})$.  Thus, using (13.3), we can find the first
$n$ terms of $\exp(P(x))$ in $\sim 9M(n)$ operations, and it is 
possible to reduce this to $\sim\frac{22}{3}M(n)$ operations by
using a fourth-order method instead of (13.3).  Compare Theorem 10.1.

\subsection*{Evaluation of $P^m$}

Suppose we want to evaluate $(P(x))^m$ for some large positive integer
$m$.  We can assume that $a_0\neq 0$, for otherwise some power of $x$
may be factored out.  Also, since $P^m = a^m_0(P/a_0)^m$, we can
assume that $a_0=1$.  By forming $P^2, P^4, P^8, \ldots \;,$ and then
the appropriate product given by the binary expansion of $m$, we can
find the first $n$ terms of $P^m$ in $O(M(n)\log_2m)$ operations.
Surprisingly, this is {\em not\/} the best possible result, at least
for large $m$.  From the identity \vspace*{-4mm}

$$ P^m = \exp(m\log(P))
\leqno{(13.4)}$$

and the above results, we can find the first $n$ terms of $P^m$
in $O(M(n))$ operations!  (If  $a_0 \neq 1$, we also need $O(\log_2 m)$
operations to evaluate $a^m_0$.)  If the methods described above
are used to compute the exponential and logarithm in (13.4), then
the number of operations is $\sim\frac{34}{3}M(n)$ as $n \ra \infty$.

\subsection*{Other operations on power series}

The method used to evaluate $\log(P(x))$ can easily be generalised to give
a method for $f(P(x))$, where $df(t)/dt$ is a function of $t$ which may be 
written in terms of square roots, reciprocals etc.  For example, with
$f(t) = \artan(t)$ we have $df/dt = 1/(1+t^2)$, so it is easy to evaluate
$\artan(P(x))$.  Using Newton's method we can evaluate the inverse
function $f^{(-1)}(P(x))$ if $f(P(x))$ can be evaluated.  Generalizations 
and applications are given in Brent and Kung $[$75$]$.\\

Some operations on formal power series do not correspond to natural
operations on multiple-precision numbers.  One example, already mentioned 
above, is formal differentiation.  Other interesting examples are 
composition and reversion.  The classical composition and reversion
algorithms, as given in Knuth $[$69$]$, are $O(n^3)$, but much faster
algorithms exist:  see Brent and Kung $[$75$]$.

\parskip 0mm

\addcontentsline{toc}{section}{References}
\section*{References}

{
\begin{description}

\item[Beeler, Gosper and Schroeppel {[72]}]
Beeler, M., Gosper, R.W., and Schroeppel, R.
``Hakmem''. Memo No.~239, M.I.T. Artificial Intelligence Lab.,
1972, 70--71.

\item[Borodin {[73]}]
Borodin, A., 
``On the number of arithmetics required to compute certain functions --
circa May 1973''.
In Complexity of Sequential and Parallel Numerical Algorithms
(ed.~by J.F.~Traub), Academic Press, New York, 1973, 149--180.

\item[Brent {[75a]}]	
Brent, R.P.,
``The complexity of multiple-precision arithmetic''.
Proc.~Seminar on Complexity of Computational Problem Solving (held
at the Australian National University, Dec.~1974), Queensland
Univ.~Press, Brisbane, 1975.

\item[Brent {[75b]}]	
Brent, R.P.,
``Fast multiple-precision evaluation of elementary functions''.
Submitted to J.~ACM.

\item[Brent and Kung {[75]}]
Brent, R.P.~and Kung, H.T.,
``Fast algorithms for reversion and composition of power series''.
To appear. (A preliminary paper appears in these Proceedings.)

\item[Gauss {[1876]}]
Gauss, C.F.,
``Carl Friedrich Gauss Werke'',
(Bd.~3), G\"{o}ttingen, 1876, 362--403.

\item[Katatsuba and Ofman {[62]}]
Karatsuba, A.~and Ofman, Y.,
``Multiplication of multidigit numbers on automata'', (in Russian).
Dokl.~Akad.~Nauk SSSR 146 (1962), 293--294.

\item [Knuth {[69]}]
Knuth, D.E.,
``The Art of Computer Programming'', (Vol.~2),
Addison Wesley, Reading, Mass., 1969, Sec.~4.7.

\item[Legendre {[11]}]
Legendre, A.M., 
``Exercices de Calcul Integral'', (Vol.~1),
Paris, 1811, 61.

\item[Melzak {[73]}]
Melzak, Z.A.,
``Companion to Concrete Mathematics'',
Wiley, New York, 1973, 68--69.

\item[Salamin {[75]}]
Salamin, E.,
``A fast algorithm for the computation of $\pi$''.
To appear in Math.\ Comp.

\item[Sch\"{o}nhage and Strassen {[71]}]
Sch\"{o}nhage, A.~and Strassen, V.,
``Schnelle Multiplikation\linebreak
grosser Zahlen''.
Computing 7 (1971), 281--292.

\item[Shanks and Wrench {[62]}]
Shanks, D.~and Wrench, J.W.,
``Calculation of $\pi$ to 100,000 decimals''.
Math.~Comp.~16 (1962), 76--99.

\end{description}
}

\newpage

\section*{Postscript (September 1999)}
\subsection*{Historical Notes and References}

This paper was retyped in \LaTeX\ with minor corrections
in September 1999.
It is available electronically in compressed postscript format from\\
{\tt http://{\lbrk}wwwmaths.anu.edu.au/{\lbrk}{\realtilde}brent/%
{\lbrk}pub/{\lbrk}pub028.html}
\medskip

Brent~[75a] is available electronically in compressed postscript format from\\
{\tt http://{\lbrk}wwwmaths.anu.edu.au/{\lbrk}{\realtilde}brent/%
{\lbrk}pub/{\lbrk}pub032.html}
\medskip

Brent~[75b]	
appeared in {\em Journal of the ACM} {\bf{23}} (1976), 242--251.
See\\
{\tt http://{\lbrk}theory.lcs.mit.edu/{\lbrk}{\realtilde}jacm/%
{\lbrk}jacm76.html\#Brent1976:242}
\medskip

The ``preliminary paper'' Brent and Kung~[75] 
appeared as ``$O((n \log n)^{3/2})$ algorithms
 for composition and reversion of power series''
 in {\em Analytic Computational Complexity}
 (edited by J.~F.~Traub),
 Academic Press, New York, 1975, 217--225.
The final paper 
appeared as
``Fast algorithms for manipulating formal power series'',
{\em Journal of the ACM} {\bf{25}} (1978), 581--595.
See\\
{\tt http://{\lbrk}theory.lcs.mit.edu/{\lbrk}{\realtilde}jacm/%
{\lbrk}jacm78.html\#BrentK1978:581}
\medskip

A generalisation to multivariate power series	
appeared as
R.~P.~Brent and H.~T.~Kung, ``Fast algorithms for composition
and reversion of multivariate power series (preliminary version)'',
in {\em Proceedings of a Conference on Theoretical Computer Science}
held at the University of Waterloo,
Dept.\ of Computer Science, University of Waterloo,
Waterloo, Ontario (August 1977),
149--158.
Abstract available electronically from\\
{\tt http://{\lbrk}wwwmaths.anu.edu.au/{\lbrk}{\realtilde}brent/%
{\lbrk}pub/{\lbrk}pub039.html}
\medskip

For more on (generalized) composition of power series, see
R.~P.~Brent and J.~F.~Traub, ``On the complexity of composition and
generalized composition of power series'',
{\em SIAM J.~Computing} {\bf{9}} (1980), 54--66.
Abstract available electronically from\\
{\tt http://{\lbrk}wwwmaths.anu.edu.au/{\lbrk}{\realtilde}brent/%
{\lbrk}pub/{\lbrk}pub050.html}
\medskip

{\em Hakmem} by Beeler, Gosper and Schroeppel~[72] is available electronically
in various formats from
{\tt http://{\lbrk}www.inwap.com/{\lbrk}pdp10/%
{\lbrk}hbaker/{\lbrk}hakmem/{\lbrk}hakmem.html}
\medskip

Salamin~[75] appeared as ``Computation of $\pi$ using arithmetic-geometric
mean'' in {\em Mathematics of Computation} {\bf{30}} (1976), 565--570.
\medskip

For much more on the arithmetic-geometric mean, see
J.~M.~Borwein and P.~B.~Borwein, {\em Pi and the AGM},
Wiley-Interscience, 1987.

\subsection*{Sharper Results}
Some of the constants can be improved.
For example, $\pi$ can be computed in $\sim 6.25M(n)\log_2n$
by the Gauss-Legendre method of Section~\ref{sec:pi},
and the constant $13$ in Theorem~9.1 can be replaced by 10.5.
For more information see the postscript to
Brent~[75a], available electronically from\\
{\tt http://{\lbrk}wwwmaths.anu.edu.au/{\lbrk}{\realtilde}brent/%
{\lbrk}pub/{\lbrk}pub032.html}

\end{document}